\documentclass{article}
\usepackage[a4paper, total={6in, 9in}]{geometry}
\usepackage[utf8]{inputenc}
\usepackage{amsmath,amsfonts}
\usepackage{booktabs} 
\usepackage{caption} 
\usepackage{subcaption} 
\usepackage{graphicx}
\usepackage{pgfplots}
\usepackage[all]{nowidow}
\usepackage[utf8]{inputenc}
\usepackage{tikz}
\usetikzlibrary{er,positioning,bayesnet}
\usepackage{multicol}
\usepackage{algpseudocode,algorithm,algorithmicx}
\usepackage{physics}
\usepackage{amssymb}
\usepackage{amsthm}
\usepackage{hyperref}
\usepackage[font=footnotesize,labelfont=bf]{caption}
\newtheorem{theorem}{Theorem}[]
\newtheorem{lemma}{Lemma}[]

\definecolor{blue}{HTML}{1F77B4}
\definecolor{orange}{HTML}{FF7F0E}
\definecolor{green}{HTML}{2CA02C}

\pgfplotsset{compat=1.14}

\begin{document}
\title{Existence of All Wilton Ripples of the Kawahara Equation}
%
%
\author{\large Ryan P. Creedon$^1$}
%
%

\date{\normalsize \today}

\maketitle              
\vspace*{-0.7cm}
{\begin{center} {\scriptsize\noindent $^1$Division of Applied Mathematics, Brown University, RI, USA, \href{mailto:ryan_creedon@brown.edu}{ryan\textunderscore creedon@brown.edu} }\end{center}}

\begin{abstract}
We investigate the existence of Wilton ripple solutions of the Kawahara equation. Without loss of generality, these are $2\pi$-periodic, traveling-wave solutions whose profiles at zero amplitude have a codimension-1 bifurcation from a linear combination of $\cos(x)$ and $\cos(Kx)$ for $K \in \mathbb{N} \setminus \{1\}$. Using a Lyapunov-Schmidt reduction, we prove the existence of these solutions for all $K$, in contrast to previous work demonstrating existence only for $K = 2$. Although the proof holds only for the Kawahara equation, many ideas introduced in the proof can be applied to more general contexts, including Wilton ripples of the gravity-capillary water wave equations.   
\textbf{Funding Acknowledgments}: This work was funded by the National Science Foundation through grant NSF-DMS-2402044. \\\\
\noindent 
\textbf{Data Availability Statement}: The Mathematica notebook \emph{CompanionToKawaharaWiltonRipples.nb} can be downloaded at \href{https://github.com/rpac5130/KawaharaWiltonRipple}{https://github.com/rpac5130/KawaharaWiltonRipple}. \\\\
\noindent \textbf{Conflict of Interest Statement}: On behalf of all authors, the corresponding author states that there is no conflict of interest.
\end{abstract}
%


%
\section{Introduction} 

We consider the Kawahara equation\footnote{Alternatively, the fifth-order KdV or super-KdV equation.} 
\begin{align}
u_t = \alpha u_{xxx} + \beta u_{5x} + \sigma \left(u^2\right)_x, \label{kEqn1}
\end{align}
where $\alpha, \beta,$ and $\sigma$ are real, nonzero parameters. First derived in \cite{kawahara1972}, the Kawahara equation models shallow water gravity-capillary waves near a critical value of the Bond number, a non-dimensional quantity that compares the forces of gravity and surface tension. In this context, the parameters $\alpha$ and $\sigma$ represent shallow-water dispersion and steepening, as in the KdV equation, while the parameter $\beta$ is a proxy for the strength of surface tension. Because of the inclusion of the surface tension term, the Kawahara equation admits several interesting and physically relevant solutions; see, for instance, \cite{alyousef2022new}-\cite{assas2009new} and \cite{mancas2019traveling,sirendaoreji2004new,ye2022all}. Of importance in this work is a class of periodic, traveling-wave solutions called the Wilton ripples.

In this work, we define Wilton ripples as a $2\pi/\kappa$-periodic, traveling-wave solution\footnote{Momentarily, we will take $\kappa = 1$ without loss of generality.} whose profile at zero amplitude has a codimension-1 bifurcation from a linear combination of $\cos(\kappa x)$ and $\cos(K \kappa x)$ for some $K \in \mathbb{N} \setminus \{1\}$. In other words, Wilton ripples represent a special class of periodic traveling waves that are parameterized by one small-amplitude parameter $a$ and have a resonant, two-mode behavior when $a = 0$, where the ratio of the modes is $1:K$. Other definitions appear in the literature: some restrict Wilton ripples to the 1:2 resonance only \cite{jones1989small,reeder1981wiltona,reeder1981wiltonb,toland1985bifurcation}, while others view them more broadly as a sheet of solutions for which the codimension of the bifurcation is greater than 1 \cite{ehrnstrom2015trimodal,kozlov2018modal,martin2013existence,wang2025existence}. Our choice of definition is consistent with that originally intended by Wilton in the context of gravity-capillary water waves \cite{akers2025ripples,trichtchenko2016instability,wilton1915lxxii} and has been used in other fields that study these solutions, including magnetohydrodynamics \cite{zakaria2004wilton}.

To motivate the mathematical behavior of Wilton ripples in the context of the Kawahara equation, we begin by searching for periodic, traveling-wave solutions. To this effect, we move to a traveling frame $x \rightarrow x - ct$ so that \eqref{kEqn1} becomes 
\begin{align}
u_t = cu_x + \alpha u_{xxx} + \beta u_{5x} + \sigma \left( \sigma^2 \right)_x.\label{kEqn_trav}
\end{align}
If $c$ is chosen to match the velocity of the traveling waves, then $u_t = 0$, and we can integrate \eqref{kEqn_trav} once in $x$, yielding
\begin{align}
cu + \alpha u_{xx} + \beta u_{4x} + \sigma u^2 = \mathcal{I}, \label{kEqn_noT_I}
\end{align}
where $\mathcal{I} \in \mathbb{R}$ is a constant of integration. By exploiting a Galilean symmetry of \eqref{kEqn1}, one can show that there exist boosts $u \rightarrow u + \xi_1$ and $c \rightarrow c + \xi_2$ of \eqref{kEqn_noT_I} such that $\mathcal{I} = 0$ without loss of generality.

What remains is to study the $2\pi/\kappa$-periodic solutions of 
\begin{align}
    cu + \alpha u_{xx} + \beta u_{4x} + \sigma u^2 = 0,
\end{align}
for all $\kappa > 0$. In fact, we can rescale $x$ and $u$ according to $x \rightarrow \kappa^{-1} x$  and $u \rightarrow \alpha k^2 \sigma^{-1} u$, respectively. Then, upon redefining the parameters $c$ and $\beta$ in terms of $\overline{c} =c/(\alpha k^2)$, and $\overline{\beta} = \kappa^2\beta/\alpha$, respectively, we can instead investigate the $2\pi$-periodic solutions of
\begin{align}
\overline{c}u + u_{xx} +  \overline{\beta} u_{4x} + u^2 = 0,\label{kEqn_final}
\end{align}
without loss of generality. In this rescaled equation, only the parameters $\overline{c}$ and $\overline{\beta}$ remain. The parameter $\overline{c}$ retains its physical interpretation as a non-dimensional velocity of the traveling wave, while $\overline{\beta}$ serves as a proxy for the Bond number. More precisely, $\overline{\beta}$ compares the relative magnitude of $\beta/\alpha$ for a fixed periodic wave. Hence, increasing $\overline{\beta}$ increases the overall effect of dispersion due to surface tension over dispersion due to gravitational effects alone. For ease of notation, the overbars on $c$ and $\beta$ will be henceforth removed. 

Because \eqref{kEqn_final} is autonomous in $x$ and satisfies the discrete symmetry $x \rightarrow -x$, we need only consider its $2\pi$-periodic, even solutions. If we further restrict these solutions to small-amplitude, \eqref{kEqn_final} effectively linearizes to
\begin{align}
cu + u_{xx} + \beta u_{4x} = 0. \label{lin_stokes}
\end{align}
At first glance, it appears that the only\footnote{It is also possible that $\cos(nx)$ solves \eqref{lin_stokes} when $c = n^2 - \beta n^4$ for any $n \in \mathbb{N}$. However, by folding $n$ into the definition of $\kappa$ before rescaling $x$, one can restrict to the $n = 1$ case without loss of generality.} real-valued solution of \eqref{lin_stokes} with the requisite parity and periodicity occurs when $c = 1 -\beta$ and is given by $u = a\cos(x)$ for any $a \in \mathbb{R}$. However, if $\beta = 1/(1+K^2)$ for any $K \in \mathbb{N} \setminus \{1\}$, then $u = a\cos(x) + b\cos(Kx)$ satisfies \eqref{lin_stokes} with $c = 1-\beta$ for all $a,b \in \mathbb{R}$. Hence, depending on the value of $\beta$, we have two classes of small-amplitude, periodic traveling-wave solutions of the Kawahara equation: 
\begin{center}
\begin{figure}[tb]
\includegraphics[width=15cm]{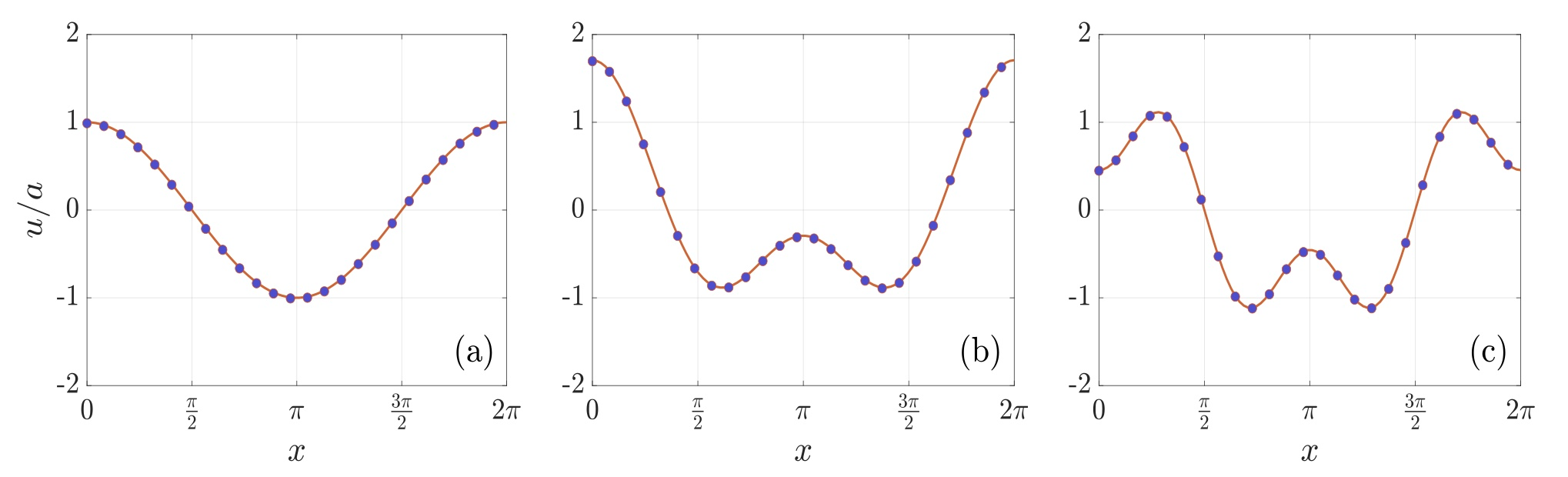}
\caption{(a) A Stokes wave of the Kawahara equation for $\beta = 1/2$, (b) A Wilton ripple of the Kawahara equation for $\beta = 1/5$, and (c) A Wilton ripple of the Kawahara equation for $\beta = 1/10$. In all plots, the amplitude parameter $a = 0.01$. The blue dots are numerically computed wave profiles using the continuation method presented in \cite{trichtchenko2018stability}. The red curves are the leading-order asymptotic behavior of the wave profiles according to Theorem 1 below. The wave profiles are normalized by $a$ for ease of comparisons.}
\end{figure}
\end{center}
 ~\vspace*{-1.21cm}
\begin{enumerate}
\item When $\beta \neq 1/(1+K^2)$, the periodic traveling-wave solutions bifurcate from a single cosine wave. For $|a|\ll 1$, these solutions take the form
\begin{align}
u = a\cos(x) + u_r(x;a) \quad \textrm{and} \quad c = 1-\beta + c_r(a),
\end{align}
where $u_r(x;a)$ is a smooth, $2\pi$-periodic function in $x$ that is real analytic in $a$ with $||u_r(x;a)||_{\textrm{L}^2} = \mathcal{O}\left(a^2\right)$ as $a \rightarrow 0$ and $c_r(a)$ is an even, real-analytic function in $a$ such that $|c_r(a)| = \mathcal{O}\left(a^2\right)$ as $a \rightarrow 0$. We call such solutions the \emph{Stokes waves} of the Kawahara equation by analogy with the corresponding solutions of the surface gravity wave equations \cite{stokes1847theory}. Their existence for sufficiently small $a$ was proven in \cite{haragus2006spectral} and their stability discussed in \cite{creedon2021high,haragus2006spectral,maspero2024full,trichtchenko2018stability}. The profile of a typical Stokes wave solution is provided in Figure 1(a). 
\item When $\beta = 1/(1+K^2)$, the periodic traveling-wave solutions bifurcate from a superposition of two co-propagating cosine waves with wavenumbers in ratio $1:K$. These are the \emph{Wilton ripples} of the Kawahara equation, according to our definition. The Wilton ripple solutions have fundamentally different expansions for small $a$ than those of the Stokes waves. Indeed, from the formal asymptotic expansions derived in \cite{haupt1988modeling,langer2024wilton}, Wilton ripple solutions of the Kawahara equation take the form
\begin{align}
u = a\Big(\cos(x) + b(a)\cos(Kx)\Big) + u_r(x;a) \quad \textrm{and} \quad c = 1-\beta + c_r(a), \label{wilton_exp}
\end{align}
where $b(a)$ is a real analytic function of $a$, $u_r(x;a)$ is a smooth, $2\pi$-periodic function of $x$ that is real analytic in $a$, and $c_r(a)$ is a real analytic function of $a$. The leading-order behaviors of $b(a)$, $u_r(x;a)$, and $c_r(a)$ as $a \rightarrow 0$ depend on $K$ according to Theorem 1 below. Profiles of typical Wilton ripple solutions are provided in Figures 1(b) and 1(c). 
\end{enumerate}
~ \hspace*{-0.4cm} \indent Existence proofs for Wilton ripples have emerged in a variety of settings, beginning with the classical gravity--capillary water wave equations. Reeder and Shinbrot~\cite{reeder1981wiltona,reeder1981wiltonb} gave the first rigorous construction of small-amplitude Wilton ripples in finite and infinite depth for the $1\!:\!2$ resonance case using a Lyapunov--Schmidt reduction. Toland and Jones~\cite{toland1985bifurcation,jones1989small} extended these considerations to small-amplitude traveling waves with codimension-1 and -2 bifurcations from flat water involving $\cos(Kx)$ and $\cos((K+1)x)$, of which the $1\!:\!2$ Wilton ripple is a special case. More recent studies have constructed multimodal resonant solutions (\emph{i.e.}, those with codimension greater than one) in various settings, including gravity--capillary waves with vorticity~\cite{ehrnstrom2015trimodal,kozlov2018modal,martin2013existence} and stratification~\cite{wang2025existence}.

Despite the considerable progress in understanding multimodal traveling waves, the existence of $1\!:\!K$ Wilton ripple solutions of the gravity--capillary water wave equations remains an open problem for general $K \in \mathbb{N} \setminus \{1\}$. The main challenge lies in tracking certain higher-order terms in the small-amplitude expansion of the Wilton ripples. Specifically, for each $K$, one must determine the order in $a$ at which the $\cos(Kx)$ mode first appears in the expansion. As $K$ increases, this order increases as well, adding to the challenge. For the full gravity--capillary equations, identifying the leading-order contribution of the $\cos(Kx)$ mode for arbitrary $K$ remains intractable. However, for simpler model equations such as the Kawahara equation, the situation is more favorable.

For the Kawahara equation and related weakly nonlinear models, the existence of $1\!:\!2$ Wilton ripples has been rigorously established via the method of majorants~\cite{akers2021wiltonB,akers2021wiltonA}. In principle, this technique can be extended to $1\!:\!K$ ripples for general $K$, but as noted in~\cite{akers2025ripples}, the complexity of the required computations grows rapidly with $K$, for essentially the same reasons as in the full water wave setting. An alternate method developed in~\cite{ehrnstrom2019bifurcation} for the gravity--capillary Whitham equation circumvents this issue by allowing the amplitudes of $\cos(x)$ and $\cos(Kx)$ to serve as free parameters, treating the remaining coefficients in the governing equation\footnote{For the gravity--capillary Whitham equation, this would be the Bond number. For the Kawahara equation, this would be $\beta$.} as dependent variables. This leads to a construction of two-dimensional \emph{sheets} of solutions for each $K$, rather than discrete branches. However, since Wilton ripples are traditionally defined for fixed values of equation parameters, it is unclear whether the classical $1\!:\!K$ branches lie embedded in these sheets. Numerical and formal asymptotic evidence from~\cite{akers2025ripples} supports this embedding when $K = 2$ and $K = 3$ for Wilton ripples of the Kawahara equation and other weakly nonlinear models, but a rigorous justification is still lacking.

Thus, the existence of $1\!:\!K$ Wilton ripples for all $K \in \mathbb{N} \setminus \{1\}$ remains unresolved, even for model equations. The objective of this work is to close this gap for the Kawahara equation. To the author's knowledge, this is the first rigorous proof of the existence of \emph{all} $1\!:\!K$ Wilton ripple solutions in a nonlinear dispersive PDE. The proof uses Lyapunov--Schmidt reduction to demonstrate analytic dependence of the Wilton ripple solutions on a small-amplitude parameter $a$. The remaining technical obstacle then is to prove that the coefficient $b(a)$ in~\eqref{wilton_exp} does not vanish identically for any $K \ge 2$. This requires expansions of the Wilton ripple solutions to arbitrarily high order in $a$, which have only recently been computed for the Kawahara equation in~\cite{langer2024wilton} using formal methods. By rigorously justifying these expansions, we prove that $b(a)$, in fact, does not vanish identically and, hence, obtain the existence of the following $1:K$ Wilton ripple solutions of the Kawahara equation:

\begin{theorem}[]
For each $K \in \mathbb{N} \setminus \{1\}$, the Kawahara equation has the following Wilton ripple solutions:
\begin{itemize}
\item If $K = 2$ (so that $\beta = 1/5$), there exist two $a$-parameter families of solutions $u_{+}(x;a)$ and $u_{-}(x;a)$ that solve \eqref{kEqn_final} for sufficiently small $a \in \mathbb{R}$. The solutions have velocity
\begin{align}
   c_{\pm}  = \frac{4}{5} - \frac{a}{\sqrt{2}}\tilde{c_r}^{\pm}(a),
    \end{align}
    respectively, where $\tilde{c_r}^{\pm}(a)$ are real analytic functions such that $\tilde{c_r}^{\pm}(0) = \pm 1. $
    The solutions take the form 
    \begin{align}
u_{\pm}(x;a) = a\Big(\cos(x) + \frac{1}{\sqrt{2}}\tilde{b}^{\pm}(a)\cos(2x)\Big) + u_r\Big(x;a,\frac{a}{\sqrt{2}}\tilde{b}^{\pm}(a),-\frac{a}{\sqrt{2}}\tilde{c_r}^{\pm}(a)\Big), 
    \end{align}
    where $\tilde{b}^{\pm}(a)$ are real analytic functions such that $\tilde{b}^{\pm}(0) = \pm 1$ and $u_r \in \textrm{H}^{\infty}_{\textrm{per},\textrm{even}}(\mathbb{T})$ is a real analytic function in its final three arguments.
    \item If $K = 3$ (so that $\beta = 1/10$), there exist three a-parameter families of solutions $u_\sigma(x;a)$ for $\sigma \in \{1,2,3\}$ that solve \eqref{kEqn_final} for sufficiently small $a \in \mathbb{R}$. The solutions have velocity
    \begin{align}
c_\sigma = \frac{9}{10} + a^2\tilde{c_r}^{(\sigma)}(a),
    \end{align}
    where $\tilde{c_r}^{(\sigma)}(a)$ are real analytic functions such that
    \begin{align}
\tilde{c_r}^{(1)}(0) \approx 4.27863 , \quad \tilde{c_r}^{(2)}(0) \approx 1.48289, \quad \textrm{and} \quad \tilde{c_r}^{(3)}(0)  \approx 0.37396 .
    \end{align}
    The solutions take the form 
    \begin{align}
u_\sigma(x;a) = a\Big(\cos(x) + \tilde{b}^{(\sigma)}(a)\cos(3x)\Big) + u_r\Big(x;a,a\tilde{b}^{(\sigma)}(a),a^2\tilde{c_r}^{(\sigma)}(a)\Big), 
    \end{align}
    where $\tilde{b}^{(\sigma)}(a)$ are real analytic functions such that 
    \begin{align}
\tilde{b}^{(1)}(0) \approx -1.78374, \quad \tilde{b}^{(2)}(0) \approx -0.54488, \quad \textrm{and} \quad \tilde{b}^{(3)}(0) \approx 0.59468,
\end{align}
    and $u_r \in \textrm{H}^{\infty}_{\textrm{per},\textrm{even}}(\mathbb{T})$ is a real analytic function in its final three arguments.
    \item If $K \geq 4$ (so that $\beta = 1/(1+K^2)$), there exists one $a$-parameter family of solutions $u(x;a)$ that solves \eqref{kEqn_final} for sufficiently small $a \in \mathbb{R}$. The solutions have velocity
    \begin{align}
c = \frac{K^2}{K^2+1} + a^2\tilde{c_r}(a),
    \end{align}
    where $\tilde{c_r}(a)$ is a real analytic function such that 
    \begin{align}
\tilde{c_r}(0) = \frac{(K^2+1)(5K^2-24)}{6K^2(K^2-4)}.
    \end{align}
    The solutions take the form
    \begin{align}
        u(x;a) = a\Big(\cos(x) + \tilde{b}(a)\cos(Kx)\Big) + u_r\Big(x;a,a\tilde{b}(a),a^2\tilde{c_r}(a)\Big),
    \end{align}
    where $\tilde{b}(a)$ is a nonzero, real analytic function such that $\tilde{b}(a) = \mathcal{O}\left(a^{K-3}\right)$ as $a \rightarrow 0$ and $u_r \in \textrm{H}^{\infty}_{\textrm{per},\textrm{even}}(\mathbb{T})$ is a real analytic function in its final three arguments.
\end{itemize}
\end{theorem}
\section{Proof of Theorem 1}
\subsection{Notation and Set-Up}
Fix $K \in \mathbb{N} \setminus \{1\}$ and let $\beta = 1/(1+K^2)$. On the interval $\mathbb{T} = [-\pi,\pi]$, define the Hilbert space
\begin{align}
\textrm{H}^4_{\textrm{per},\textrm{even}}(\mathbb{T}) := \left\{f(x) = \sum_{k = 0}^{\infty} \hat{f}_k\cos(kx) \quad \textrm{s.t.} \quad \sum_{k=0}^{\infty} (1+k^2)^4\left|\widehat{f}_k\right|^2 < +\infty \right\}, \label{H4}
\end{align}
where 
\begin{align}
\widehat{f}_k = \mathcal{F}_k[f(x)] = \begin{cases} \frac{1}{\pi}\int_{-\pi}^{\pi} f(x)\cos(kx) dx & k > 0 \\
 \frac{1}{2\pi}\int_{-\pi}^{\pi} f(x)dx & k = 0
\end{cases},
\end{align}
are the Fourier cosine coefficients of $f$. Define the nonlinear operator
$F:\textrm{H}^4_{\textrm{per},\textrm{even}}(\mathbb{T}) \times \mathbb{R} \rightarrow \textrm{L}^2_{\textrm{per},\textrm{even}}(\mathbb{T})$ according to
\begin{align}
F(u,c) = (c+ \mathcal{L})u + u^2, \quad \textrm{where} \quad \mathcal{L} := \partial_x^2 + \beta\partial_x^4.\label{Fdefn}
\end{align}
Then, \eqref{kEqn_final} becomes the operator equation
\begin{align}
F(u,c) = 0. \label{opEqn}
\end{align}
Over the given domain $\textrm{H}^4_{\textrm{per},\textrm{even}}(\mathbb{T})$, the linear operator $\mathcal{L}$ is bounded for all $c \in \mathbb{R}$ and is self-adjoint with respect to the standard $\textrm{L}^2(\mathbb{T})$ inner-product:
\begin{align}
\left<v,w \right> = \int_{-\pi}^{\pi}v(x)w(x) dx. \label{innerProd}
\end{align}
As shown in the Introduction, $c_0+\mathcal{L}$ has a two-dimensional nullspace spanned by $\cos(x)$ and $\cos(K x)$ when $c_0 := 1-\beta$. Owing to self-adjointness, the cokernel of $c_0+\mathcal{L}$ is also spanned by $\cos(x)$ and $\cos(Kx)$, which guarantees that $c_0+\mathcal{L}$ is a Fredholm operator of index zero. 

In what follows, we decompose $c$ according to 
\begin{align}
c = c_0 + c_r, \label{cDecomp}
\end{align}
where $c_r \in \mathbb{R}$ is conceived as a small correction to $c_0$. We also decompose $u \in \textrm{H}^4_{\textrm{per,even}}(\mathbb{T})$ by
\begin{align}
u = u_0 + u_r, \label{uDecomp}
\end{align}
where $u_0 \in \textrm{Null}\{c_0+\mathcal{L} \}$ and $u_r \in \textrm{Null}\{c_0+\mathcal{L}_0 \}^\perp$. (Here, orthogonality is taken with respect to the inner-product \eqref{innerProd}.) By design,
\begin{align}
u_0(x) = a\cos(x) + b \cos(Kx) \label{u0}
\end{align}
for some $a,b \in \mathbb{R}$, while $u_r$ is treated as a small correction to $u_0$ involving Fourier cosines modes orthogonal to $\cos(x)$ and $\cos(Kx)$. 

Substituting \eqref{cDecomp} and \eqref{uDecomp} into \eqref{Fdefn} yields a new operator equation
\begin{align}
G(u_r,a,b,c_r)= 0, \label{Gdefn}
\end{align}
where $G : \textrm{Null}\{c_0+\mathcal{L} \}^\perp  \times \mathbb{R}^3 \rightarrow \textrm{L}^2_{\textrm{per},\textrm{even}}(\mathbb{T}) $ such that
\begin{align}
G(u_r,a,b,c_r) &:= F\left(a\cos(x)+b\cos(Kx)+u_r,c_0+c_r\right) \nonumber \\ &~=  (c_0+\mathcal{L})u_r + c_r\left(a\cos(x)+b\cos(Kx)+u_r\right) + \left(a\cos(x)+b\cos(Kx)+u_r\right)^2.
\end{align}
Equation \eqref{Gdefn} will be the starting point of our Lyapunov-Schmidt reduction. 

An intermediate step of the reduction will show the existence of a function $u_r \in \textrm{Null}\{c_0+\mathcal{L}_0\}^\perp$ that has real analytic dependence on parameters $a$, $b$, and $c_r$. We will require a Taylor series expansion of this $u_r$ around $(a,b,c_r) = (0,0,0)$, and norms for the remainders of this series will be expressed in big-oh notation. For clarity, given any $f \in X$, where $X$ is a Banach space\footnote{Typically, $\textrm{H}^4_{\textrm{per,even}}(\mathbb{T})$ or $\mathbb{R}^3$.}, and given any $n \in \mathbb{R}^+$, we say
\begin{align}
||f||_{X} = \mathcal{O}\left( \left(|a|+|b|+|c_r|\right)^n\right) \quad \textrm{as} \quad a,b, c_r \rightarrow 0,
\end{align}
provided there exist $C = C(n) > 0$ and $\delta  = \delta(n) > 0$ such that 

\begin{align}
||f||_{X} \leq C \left(|a|+|b|+|c_r|\right)^n \quad \textrm{for all} \quad |a|+|b|+|c_r| < \delta. 
\end{align}
 With the necessary notation established, we proceed with Lyapunov-Schmidt reduction on \eqref{Gdefn}.
\subsection{Solving the Auxiliary Equation}
We introduce the operators $\mathcal{P}$ and $\mathcal{Q}$ according to
\begin{subequations}
\begin{align}
\mathcal{P}f &= \widehat{f}_1\cos(x) + \widehat{f}_K\cos(Kx), \\
\mathcal{Q}f &= f - \mathcal{P}f,\end{align}
\end{subequations}
which are valid for all $f \in \textrm{L}^2_{\textrm{per,even}}(\mathbb{T})$. Here, the operator $\mathcal{P}$ projects any function in $\textrm{L}^2_{\textrm{per,even}}(\mathbb{T})$ onto the cokernel of $c_0 + \mathcal{L}$, while $\mathcal{Q}$ projects onto the range of $c_0+\mathcal{L}$.

Applying $\mathcal{Q}$ to \eqref{Gdefn} yields a new operator equation
\begin{align}
\mathcal{A}(u_r,a,b,c_r) = 0, \label{Aeqn}
\end{align}
where
\begin{align}
\mathcal{A}(u_r,a,b,c_r) :=\mathcal{Q}G(u_r,a,b,c_r) = c_ru_r + \mathcal{Q}(c_0+\mathcal{L})u_r + \mathcal{Q}\left(a\cos(x)+b\cos(Kx) + u_r\right)^2. \label{Adefn}
\end{align}
We call \eqref{Aeqn} the \textit{auxiliary equation}. From \eqref{Adefn}, it follows that 
\begin{align}
\mathcal{A}(0,0,0,0) = 0.
\end{align}
Moreover, since $F$ in \eqref{Fdefn} is real analytic\footnote{Indeed, $\mathcal{L}$ is a bounded linear operator from $\textrm{H}_{\textrm{per,even}}^4(\mathbb{T})$ to $\textrm{L}_{\textrm{per,even}}^2(\mathbb{T})$, and $u^2$ is an analytic operator from $\textrm{H}_{\textrm{per,even}}^4(\mathbb{T})$ to $\textrm{L}_{\textrm{per,even}}^2(\mathbb{T})$ because $\textrm{H}_{\textrm{per,even}}^4(\mathbb{T})$ is a Banach algebra.} with respect to $u$ and $c$ and since the projections $\mathcal{P}$ and $\mathcal{Q}$ are bounded linear operators, $\mathcal{A}$ is real analytic with respect to $u_r$, $a$, $b$, and $c_r$. 

Using analyticity, a straightforward calculation shows that 
\begin{align}
\partial_{u_r} \mathcal{A}(0,0,0) = \mathcal{Q}(c_0+\mathcal{L}),
\end{align}
where $\partial_{u_r}$ denotes the Fr\'{e}chet derivative with respect to $u_r$. The resulting operator $\mathcal{Q}(c_0+\mathcal{L}) : \textrm{Null}\{c_0+\mathcal{L}\}^\perp \rightarrow \textrm{Range}\{c_0+\mathcal{L}\}$ has bounded inverse
\begin{align}
\left(\mathcal{Q}(c_0+\mathcal{L})\right)^{-1}(f(x)) = \sum_{k \in \mathbb{N}_0 \setminus \{1,K\}} \left(\frac{\widehat{f}_k}{c_0-k^2+\beta k^4} \right) \cos(kx).
\end{align}
Thus, we have $\mathcal{A}: \textrm{Null}\left(c_0+\mathcal{L}\right)^\perp \times \mathbb{R}^3 \rightarrow \textrm{Null}\left(c_0+\mathcal{L}\right)^\perp$ such that
\begin{enumerate}
\item[(\emph{i})] $\textrm{Null}\left(c_0+\mathcal{L}\right)^\perp$ with the $\textrm{H}^4$ norm and $\mathbb{R}^3$ with the $\ell^1$ norm (for convenience) are Banach spaces, 
\item[(\emph{ii})]  $\mathcal{A}$ is real analytic with respect to both $u_r$ and $(a,b,c_r)^T$,
\item[(\emph{iii})]  $\mathcal{A}$ satisfies $\mathcal{A}(0,0,0,0) = 0$, and 
\item[(\emph{iv})] $\partial_{u_r}\mathcal{A}(0,0,0,0)$ is boundedly invertible.
\end{enumerate}
 Hence, we can apply the analytic version of the Implicit Function Theorem in Banach spaces \cite[Theorem 2.3]{chow2012methods} to deduce the existence and uniqueness of 
\begin{align}
u_r = u_r(x;a,b,c_r), \label{ur_soln}
\end{align}
such that, for sufficiently small $a$, $b$, and $c_r$, we have 
\begin{enumerate}
    \item[(\emph{i})] $u_r(x;a,b,c_r) \in \textrm{Null}\{c_0+\mathcal{L}\}^\perp \subset \textrm{H}^4_{\textrm{per},\textrm{even}}(\mathbb{T})$,
    \item[(\emph{ii})] $\mathcal{A}(u_r(x;a,b,c_r),a,b,c_r) = 0$,
    \item[(\emph{iii})]  $u_r(x;0,0,0) = 0$, and
    \item[(\emph{iv})] $u_r(x;a,b,c_r)$ is real analytic with respect to $a$, $b$, and $c_r$.
\end{enumerate}

As \eqref{ur_soln} is a real analytic map from $\mathbb{R}^3$ to $\textrm{Null}\left(c_0+\mathcal{L}\right)^\perp$, we can use Taylor's Theorem in the Banach space setting \cite[Theorem 5.4]{coleman2012calculus}  to expand $u_r$ as a multi-variable power series in $a$, $b$, and $c_r$ with values in $\textrm{Null}\left(c_0+\mathcal{L}\right)^\perp$. In particular, we have
\begin{align}
u_r(x;a,b,c_r) &= u_r^{(1,0,0)}(x)a +  u_r^{(0,1,0)}(x)b +  u_r^{(0,0,1)}(x)c_r + u_r^{(2,0,0)}(x)a^2 + u_r^{(1,1,0)}(x)ab \nonumber\\&\quad+ u_r^{(0,2,0)}(x) b^2 + u_r^{(0,1,1)}(x)bc_r + u_r^{(0,0,2)}(x) c_r^2 + u_r^{(1,0,1)}(x)ac_r + R_3(x;a,b,c_r), \label{ur_exp}
\end{align}
where $u_r^{(i,j,k)}(x)$ are explicitly computed in the Appendix and $R_3(x;a,b,c_r) \in \textrm{Null}\{c_0+\mathcal{L}\}^\perp \subset \textrm{H}^4_{\textrm{per},\textrm{even}}(\mathbb{T})$ such that
\begin{align}\Big|\Big|R_3(x;a,b,c_r)\Big|\Big|_{\textrm{H}^4(\mathbb{T})} = \mathcal{O}\left((|a|+|b|+|c_r|)^3\right) \quad \textrm{as} \quad a,b,c_r \rightarrow 0.
\end{align}
\subsection{Solving the Bifurcation Equation}
Applying the projection $\mathcal{P}$ to \eqref{Gdefn} yields a new operator equation
\begin{align}
\mathcal{B}(u_r,a,b,c_r) = 0, \label{Beqn}
\end{align}
where $\mathcal{B}: \textrm{Null}\{c_0+\mathcal{L}\}^\perp 
 \cross \mathbb{R}^3 \rightarrow \textrm{Coker}\{c_0 + \mathcal{L} \}$ such that
\begin{align}
\mathcal{B}(u_r,a,b,c_r) &:= \mathcal{P}G(u_r,a,b,c_r) \nonumber \\ &~= c_r\left(a\cos(x) + b\cos(Kx) \right) + \mathcal{P}\left(a\cos(x) + b\cos(Kx) + u_r \right)^2.
\end{align}
Substituting $u_r = u_r(x;a,b,c_r)$ into \eqref{Beqn} and separating terms proportional to $\cos(x)$ and $\cos(Kx)$, we arrive at a finite-dimensional operator equation 
\begin{align}
\tilde{\mathcal{B}}(a,b,c_r) = 0, \label{tildeBeqn}
\end{align}
where $\tilde{B}: \mathbb{R}^3 \rightarrow \mathbb{R}^2$ such that
\begin{align}
\tilde{\mathcal{B}}(a,b,c_r) = \begin{pmatrix} \mathcal{F}_1\left[\left(a\cos(x)+b\cos(Kx) + u_r(x;a,b,c_r) \right)^2\right] + ac_r \\
\mathcal{F}_K\left[\left(a\cos(x)+b\cos(Kx) + u_r(x;a,b,c_r) \right)^2\right]+ bc_r \end{pmatrix}.
\end{align}
We call \eqref{tildeBeqn} the \emph{bifurcation equation}. 

Owing to the analyticity of $u_r(x;a,b,c_r)$ with respect to $(a,b,c_r)^T$, we expand $\tilde{\mathcal{B}}$ for sufficiently small $a$, $b$, and $c_r$. To this end, we let $\omega \in \{1,K\}$ and consider
\begin{align}
\mathcal{F}_{\omega}\left[\left(a\cos(x)+b\cos(Kx) + u_r(x;a,b,c_r) \right)^2 \right] = \mathcal{T}_{\omega,1} + \mathcal{T}_{\omega,2} + \mathcal{T}_{\omega,3},
\end{align}
where
\begin{subequations}
\begin{align}
\mathcal{T}_{\omega,1} &= \mathcal{F}_{\omega}\left[\left(a\cos(x)+b\cos(Kx) \right)^2 \right], \label{T1}\\
\mathcal{T}_{\omega,2} &= \mathcal{F}_{\omega}\left[2\left(a\cos(x)+b\cos(Kx) \right)u_r(x;a,b,c_r) \right], \label{T2} \\
\mathcal{T}_{\omega,3} &= \mathcal{F}_{\omega}\left[u_r(x;a,b,c_r)^2 \right].
\end{align}
\end{subequations}
For sufficiently small $a$, $b$, and $c_r$, we have 
\begin{align}
\left|\mathcal{T}_{\omega,3} \right| \leq ||u_r(x;a,b,c_r)||^2_{\textrm{L}^2(\mathbb{T})} \leq ||u_r(x;a,b,c_r)||^2_{\textrm{H}^4(\mathbb{T})}.
\end{align}
Moreover, because the linear terms in the Taylor series of $u_r(x;a,b,c_r)$ vanish for all $K \in \mathbb{N} \setminus \{1\}$ by \eqref{urlin1} and \eqref{urlin2} in the appendix, we must have 
\begin{align}
||u_r(x;a,b,c_r)||_{\textrm{H}^4(\mathbb{T})} \leq C\left(|a|+|b|+|c_r| \right)^2,
\end{align}
for some constant $C > 0$ by Taylor's Theorem\footnote{Note that this constant must be independent of $a$, $b$, and $c_r$, but it may depend on $K$.}. Thus, 
\begin{align}
\mathcal{T}_{\omega,3} = \mathcal{O}\left(|a|+|b|+|c_r| \right)^4 \quad \textrm{as} \quad a,b,c_r \rightarrow 0, \label{T3_exp_final}
\end{align}
regardless of $\omega$. In contrast, $\mathcal{T}_{\omega,1}$ and $\mathcal{T}_{\omega,2}$ contribute terms lower than fourth degree in $a$, $b$, and $c_r$ when expanded.  Ultimately, the differences in the leading-order behavior of $\mathcal{T}_{\omega,1}$ and $\mathcal{T}_{\omega,2}$ determine the number of distinct branches of Wilton ripple solutions when $K = 2$, $K = 3$, and $K \geq 4$. Indeed, we will find two branches of real-valued solutions to the bifurcation equation  when $K = 2$ and, hence, two branches of Wilton ripple solutions in this case, while three branches are found when $K = 3$ and only one branch is found when $K \geq 4$. 
\subsubsection{Case 1: $K = 2$} 
A direct calculation of \eqref{T1} when $K = 2$ yields
\begin{align}
\mathcal{T}_{\omega,1} = \begin{cases} ab & \omega = 1 \\ \frac12 a^2 & \omega = 2\end{cases}. \label{2T1_exp_final}
\end{align}
Turning to $\mathcal{T}_{\omega,2}$, we substitute the expansion \eqref{ur_exp} into \eqref{T2}, giving us
\begin{align}
\mathcal{T}_{\omega,2} &= \mathcal{F}_{\omega}\left[2\big(a\cos(x)+b\cos(2x)\big)\big(u^{(2,0,0)}(x)a^2 + u^{(1,1,0)}(x)ab + u^{(0,2,0)}(x)b^2  \big) \right] +\mathcal{R}_{\omega,2}, \label{2T2_exp}
\end{align}
where $u^{(i,j,k)}(x)$ are given in \eqref{2u200}-\eqref{2u020} and
\begin{align}
\mathcal{R}_{\omega,2} =  \mathcal{F}_{\omega}\left[2\big(a\cos(x)+b\cos(2x)\big)R_3(x;a,b,c_r) \right].
\end{align}
Using a combination of the Cauchy-Schwarz Inequality and Taylor's Theorem, we find
\begin{align}
|\mathcal{R}_{\omega,2}| &\leq 2\big(|a|+|b| \big) ||R_3(x;a,b,c_r)||_{\textrm{L}^1(\mathbb{T})}, \nonumber \\ &\leq 2\sqrt{2\pi}(|a|+|b|)||R_3(x;a,b,c_r)||_{\textrm{L}^2(\mathbb{T})},\nonumber \\
&\leq  2\sqrt{2\pi}(|a|+|b|)||R_3(x;a,b,c_r)||_{\textrm{H}^4(\mathbb{T})} ,\nonumber \\
&\leq 2\sqrt{2\pi}C(|a|+|b|+|c_r|)^3(|a|+|b|), \nonumber \\
&\leq 2\sqrt{2\pi}C(|a|+|b|+|c_r|)^4,
\end{align}
which is valid for sufficiently small $a$, $b$, and $c_r$. It follows that
\begin{align}
\mathcal{R}_{\omega,2} = \mathcal{O}\left(|a|+|b|+|c_r| \right)^4 \quad \textrm{as} \quad a,b,c_r \rightarrow 0.
\end{align}
The remaining term in \eqref{2T2_exp} can be calculated explicitly with the aid of Mathematica\footnote{See the Mathematica companion file \emph{CompanionToKawaharaWiltonRipples.nb} for details.}. We find
\begin{align}
\mathcal{T}_{\omega,2} = \begin{cases} -\frac54 a^3 - \frac{11}{8}ab^2 + \mathcal{O}\left(|a|+|b|+|c_r| \right)^4 & \omega = 1 \\
-\frac{11}{8}a^2b - \frac{91}{72}b^3 + \mathcal{O}\left(|a|+|b|+|c_r| \right)^4 & \omega = 2\end{cases}, \label{2T2_exp_final}
\end{align}
where we have used $\beta = 1/(1+K^2) = 1/5$ and $c_0 = 1-\beta = 4/5$ to simplify. Substituting \eqref{T3_exp_final}, \eqref{2T1_exp_final}, and \eqref{2T2_exp_final} into \eqref{tildeBeqn} yields the following expansion of the bifurcation equation:
\begin{align}
\tilde{\mathcal{B}}(a,b,c_r) &= \begin{pmatrix} ab + ac_r - \frac54 a^3 - \frac{11}{8}ab^2 + \mathcal{O}\left(|a|+|b|+|c_r|\right)^4 \\ \frac12 a^2 +bc_r - \frac{11}{8}a^2b - \frac{91}{72}b^3 +  \mathcal{O}\left(|a|+|b|+|c_r|\right)^4 \end{pmatrix} = 0. \label{BtildeExp}
\end{align}
Comparison of the dominant terms in \eqref{BtildeExp} suggests the following rescaling:
\begin{align}
b = \frac{a}{\sqrt{2}}\tilde{b}, \quad c_r = -\frac{a}{\sqrt{2}}\tilde{c_r}, \quad \textrm{and} \quad \tilde{\tilde{\mathcal{B}}}(a,\tilde{b},\tilde{c_r}) = \frac{1}{a^2}\tilde{\mathcal{B}}\left(a,\frac{a}{\sqrt{2}}\tilde{b},-\frac{a}{\sqrt{2}}\tilde{c_r}\right).
\end{align}
Then, \eqref{BtildeExp} becomes
\begin{align}
\tilde{\tilde{B}}(a,\tilde{b},\tilde{c_r}) = \begin{pmatrix} \frac{1}{\sqrt{2}}(\tilde{b} -\tilde{c_r})- \frac54 a - \frac{11}{16}a\tilde{b}^2 + \mathcal{O}\left(|a|^2\left(1 + \frac{1}{\sqrt{2}}|\tilde{b}| + \frac{1}{\sqrt{2}}|\tilde{c_r}| \right)^4\right) \\ \frac12(1-\tilde{b}\tilde{c_r}) - \frac{11}{8\sqrt{2}}a\tilde{b} - \frac{91}{144\sqrt{2}}a\tilde{b}^3 + \mathcal{O}\left(|a|^2\left(1 + \frac{1}{\sqrt{2}}|\tilde{b}| + \frac{1}{\sqrt{2}}|\tilde{c_r}| \right)^4\right) \end{pmatrix} = 0. \label{BtildetildeExp}
\end{align}
From \eqref{BtildetildeExp}, we have that $\tilde{\tilde{\mathcal{B}}}(0,\pm1,\pm1) = 0$.  The corresponding Jacobian evaluated at $(0,\pm1,\pm1)$ is
\begin{align}
 \frac{\partial \tilde{\tilde{\mathcal{B}}}}{\partial(\tilde{b},\tilde{c_r})}(0,\pm1,\pm1) = \begin{pmatrix} \frac{1}{\sqrt{2}} & -\frac{1}{\sqrt{2}} \\ \mp 1 & \mp 1 \end{pmatrix},
\end{align}
which is invertible in both cases. The Analytic Implicit Function Theorem \cite[Theorem 2.3]{chow2012methods} guarantees the existence and uniqueness of 
\begin{align}
\tilde{b}^{\pm} = \tilde{b}^{\pm}(a) \quad \textrm{and} \quad
\tilde{c_r}^{\pm} = \tilde{c_r}^{\pm}(a),\end{align}
such that, for sufficiently small $a$, we have (\emph{i}) $\tilde{\tilde{\mathcal{B}}}(a,\tilde{b}^{\pm}(a),\tilde{c_r}^{\pm}(a)) = 0$, (\emph{ii}) $\tilde{b}^{\pm}(0) = \pm 1$ and $\tilde{c_r}^{\pm}(0) = \pm 1$, and (\emph{iii})  $\tilde{b}^{\pm}(a)$ and $\tilde{c_r}^{\pm}(a)$ are real analytic in $a$.
    
    Having solved both the auxiliary and bifurcation equations when $K = 2$, we conclude the existence of sufficiently small-amplitude Wilton ripple solutions of \eqref{kEqn_final} in this case. In particular, if $K = 2$ and $\beta = 1/(1+K^2) = 1/5$, there exist two unique $a$-parameter family of solutions to \eqref{kEqn_final} for sufficiently small $a\in \mathbb{R}$. These solutions travel at velocity \begin{align}
   c(a) = c_0 - \frac{a}{\sqrt{2}}\tilde{c_r}^{\pm}(a) = \frac{4}{5} - \frac{a}{\sqrt{2}}\tilde{c_r}^{\pm}(a),
    \end{align}
    respectively, and take the form
    \begin{align}
u_{\pm}(x;a) = a\Big(\cos(x) + \frac{1}{\sqrt{2}}\tilde{b}^{\pm}(a)\cos(2x)\Big) + u_r\Big(x;a,\frac{a}{\sqrt{2}}\tilde{b}^{\pm}(a),-\frac{a}{\sqrt{2}}\tilde{c_r}^{\pm}(a)\Big), \label{usolnK2}
    \end{align}
    where $u_r \in \textrm{H}^4_{\textrm{per},\textrm{even}}(\mathbb{T})$ is real analytic in its parameters with Taylor expansion \eqref{ur_exp} and coefficients \eqref{2u200}-\eqref{2u020}. The functions $\tilde{b}^{\pm}(a)$ and $\tilde{c_r}^{\pm}(a)$ are real analytic in $a$ such that $\tilde{b}^\pm(0) = \pm 1$ and $\tilde{c_r}(0) = \pm 1$. Consequently, $u_{\pm}(x;a)$ are real analytic in $a$ with values in $\textrm{H}^4_{\textrm{per},\textrm{even}}(\mathbb{T})$. In fact, by rewriting \eqref{kEqn_final} as
    \begin{align}
    \partial_x^4 u_{\pm}(x;a) = -\beta^{-1}\left(cu_{\pm}(x;a) + \partial_x^2u_{\pm}(x;a) + u_{\pm}(x;a)^2 \right),
    \end{align}
    we infer $\partial_x^4u_{\pm}(x;a) \in \textrm{H}^{2}_{\textrm{per},\textrm{even}}(\mathbb{T})$ and, hence, $u_{\pm}(x;a) \in \textrm{H}^{6}_{\textrm{per},\textrm{even}}(\mathbb{T})$. Repeating this argument indefinitely, one concludes $u_{\pm}(x;a) \in \textrm{H}^{\infty}_{\textrm{per},\textrm{even}}(\mathbb{T})$, completing the proof of Theorem 1 when $K = 2$. 
    \subsubsection{Case 2: $K = 3$}
    When $K = 3$, a direct calculation shows
    \begin{align}
\mathcal{T}_{\omega,1} = 0,
    \end{align}
    regardless of $\omega$. Turning to $\mathcal{T}_{\omega,2}$, we substitute \eqref{u2_k}-\eqref{u02_k} into \eqref{2T2_exp} with $K = 3$ to obtain the following expansion:
    \begin{align}
\mathcal{T}_{\omega,2} = \begin{cases} -\frac79 a^3 + a^2b - \frac{34}{63} ab^2 + \mathcal{O}\left(|a|+|b|+|c_r| \right)^4 & \omega = 1 \\ \frac13 a^3 - \frac{34}{63}a^2b - \frac{211}{189} b^3 +  \mathcal{O}\left(|a|+|b|+|c_r| \right)^4 & \omega = 3\end{cases}, \label{3T2_exp_final}
    \end{align}
    where we have used $\beta = 1/(1+K^2) = 1/10$ and $c_0 = 1-\beta = 9/10$ to simplify the coefficients above. The remainder terms are estimated exactly as in the previous case.

    Given \eqref{3T2_exp_final}, the bifurcation equation expands as
    \begin{align}
\tilde{\mathcal{B}}(a,b,c_r) &= \begin{pmatrix} ac_r -\frac79 a^3 + a^2b - \frac{34}{63} ab^2 + \mathcal{O}\left(|a|+|b|+|c_r| \right)^4 \\ bc_r +  \frac13 a^3 - \frac{34}{63}a^2b - \frac{211}{189} b^3 +  \mathcal{O}\left(|a|+|b|+|c_r| \right)^4 \end{pmatrix} = 0. \label{BtildeExp3}
\end{align}
Balancing the quadratic and cubic terms above suggests the following rescaling:
\begin{align}
b = a\tilde{b}, \quad c_r = a^2\tilde{c_r}, \quad \textrm{and} \quad \tilde{\tilde{\mathcal{B}}}(a,\tilde{b},\tilde{c_r}) = \frac{1}{a^3}\tilde{B}(a,a\tilde{b},a^2\tilde{c_r}), \label{rescale3}
\end{align}
leading us to a new bifurcation equation
\begin{align}
\tilde{\tilde{\mathcal{B}}}(a,\tilde{b},\tilde{c_r}) = \begin{pmatrix} -\frac79 + \tilde{c_r} + \tilde{b} - \frac{34}{63}\tilde{b}^2 + \mathcal{O}\left(|a|\left(1+|\tilde{b}|+|a|^2|\tilde{c_r}| \right)^4 \right) \\ \frac13 - \frac{34}{63}\tilde{b} + \tilde{b}\tilde{c_r} - \frac{211}{189}\tilde{b}^3 +  \mathcal{O}\left(|a|\left(1+|\tilde{b}|+|a|^2|\tilde{c_r}| \right)^4 \right)\end{pmatrix} = 0.
\end{align}
When $a = 0$, we obtain a solution to the above provided
\begin{align}
\tilde{c_r} = \frac79 -\tilde{b} + \frac{34}{63}\tilde{b}^2 \quad \textrm{and} \quad \tilde{b}^3+\tilde{b}^2 - \frac{5}{21}\tilde{b} - \frac13 = 0. \label{crbeqns}
\end{align}
Since the discriminant of the $\tilde{b}$  equation is positive, we anticipate three distinct, real solutions. Arranged from least to greatest, we denote these solutions by $\tilde{b}^{(\sigma)}_0$ for $\sigma = 1,~ 2, ~\textrm{and} ~3$, respectively.  Numerical methods reveal
\begin{align}
\tilde{b}^{(1)}_0 \approx -1.78374, \quad \tilde{b}^{(2)}_0 \approx -0.54488, \quad \textrm{and} \quad \tilde{b}^{(3)}_0 \approx 0.59468. \label{btildenum}
\end{align}
For each $\sigma$, denote $\tilde{c_{r,0}}^{(\sigma)}$ by the evaluation of \eqref{crbeqns} at $\tilde{b}_0^{(\sigma)}$. Using the numerical values \eqref{btildenum}, we find
\begin{align}
\tilde{c_{r,0}}^{(1)} \approx 4.27863 ,\quad \tilde{c_{r,0}}^{(2)} \approx 1.48289, \quad \textrm{and} \quad \tilde{c_{r,0}}^{(3)} \approx 0.37396.
\end{align}
By construction, it follows necessarily that 
\begin{align}
\tilde{\tilde{\mathcal{B}}}\left(0,\tilde{b}_0^{(\sigma)},\tilde{c_{r,0}}^{(\sigma)}\right) = 0,
\end{align}
for each $\sigma \in \{1,2,3\}$. The corresponding Jacobians are
\begin{align}
\frac{\partial \tilde{\tilde{\mathcal{B}}}}{\partial(\tilde{b},\tilde{c_{r}})}(0,\tilde{b}_0^{(\sigma)},\tilde{c_{r,0}}^{(\sigma)}) = \begin{pmatrix} 1 & 1 \\ -\frac{34}{63} & \tilde{b}_0^{(\sigma)} \end{pmatrix},
\end{align}
which are invertible since $\tilde{b}_0^{(\sigma)} \neq -34/63$ for all $\sigma$ by contradiction. The Analytic Implicit Function Theorem then gives the existence and uniqueness of 
\begin{align}
\tilde{b}^{(\sigma)} = \tilde{b}^{(\sigma)}(a) \quad \textrm{and} \quad \tilde{c_{r}}^{(\sigma)} = \tilde{c_r}^{(\sigma)}(a),
\end{align}
 for each $\sigma \in \{1,2,3\}$ and sufficiently small $a$ such that (\emph{i}) $\tilde{\tilde{\mathcal{B}}}(a,\tilde{b}^{(\sigma)}(a),\tilde{c_r}^{(\sigma)}(a)) = 0$, (\emph{ii}) $\tilde{b}^{(\sigma)}(0) = \tilde{b}_0^{(\sigma)}$ and $\tilde{c_r}^{(\sigma)}(0) = \tilde{c_{r,0}}^{(\sigma)}$, and (\emph{iii}) $\tilde{b}^{(\sigma)}(a)$ and $\tilde{c_r}^{(\sigma)}(a)$ are real analytic in $a$.

     We conclude the existence of sufficiently small-amplitude Wilton ripples of \eqref{kEqn_final} when $K = 3$, \emph{i.e.}, $\beta = 1/10$. In this case, there exist three unique $a$-parameter family of solutions to \eqref{kEqn_final} for sufficiently small $a \in \mathbb{R}$. These solutions travel at velocity
     \begin{align}
c(a) = c_0 + a^2\tilde{c_{r}}^{(\sigma)}(a) = \frac{9}{10} + a^2\tilde{c_{r}}^{(\sigma)}(a) , 
     \end{align}
    for $\sigma \in \{1,2,3\}$, and their corresponding profiles are
    \begin{align}
u_{\sigma}(x;a) = a\Big(\cos(x) + \tilde{b}^{(\sigma)}(a)\cos(3x)\Big) + u_r\Big(x;a,a\tilde{b}^{(\sigma)}(a),a^2\tilde{c_r}^{(\sigma)}(a)\Big),
    \end{align}
     where $u_r \in \textrm{H}^4_{\textrm{per},\textrm{even}}(\mathbb{T})$ is real analytic in its parameters with Taylor expansion \eqref{ur_exp} and coefficients \eqref{u2_k}-\eqref{u02_k} evaluated at $K = 3$. The functions $\tilde{b}^{(\sigma)}(a)$ and $\tilde{c_r}^{(\sigma)}(a)$ are real analytic for sufficiently small $a$ and satisfy
     \begin{align}
\tilde{b}^{(\sigma)}(0) = \tilde{b}_0^{(\sigma)} \quad \textrm{and} \quad  \tilde{c_r}^{(\sigma)}(0) = \tilde{c_{r,0}}^{(\sigma)}.
     \end{align}
     By the same bootstrap argument as before, we conclude $u_{\sigma}(x;a)$ are real analytic functions of $a$ with values in $\textrm{H}^{\infty}_{\textrm{per},\textrm{even}}(\mathbb{T})$ for each $\sigma$, completing the proof of Theorem 1 when $K = 3$.
\subsubsection{Case 3: $K \geq 4$}
  As in the $K = 3$ case, we find
    \begin{align}
\mathcal{T}_{\omega,1} = 0,
    \end{align}
    regardless of $\omega$. Substituting \eqref{u2_k}-\eqref{u02_k} into \eqref{2T2_exp} yields the following expansion for $K \geq 4$:
    \begin{align}
\mathcal{T}_{\omega,2} = \begin{cases} v_{3,0,0}(K)a^3 + v_{1,2,0}(K)ab^2 + \mathcal{O}\left(|a|+|b|+|c_r| \right)^4 & \omega = 1 \\ w_{2,1,0}(K)a^2b +w_{0,3,0}(K)b^3  + \mathcal{O}\left(|a|+|b|+|c_r| \right)^4 & \omega = K\end{cases}, \label{KT2_exp_final}
    \end{align}
    where 
    \begin{subequations}
    \begin{align}
v_{3,0,0}(K) &= -\frac{(K^2+1)(5K^2-24)}{6K^2(K^2-4)}, \\ 
v_{1,2,0}(K) = w_{2,1,0}(K) &= -\frac{(K^2+1)(4K^4-27K^2+4)}{K^2(K^2-4)(4K^2-1)}, \\
w_{0,3,0}(K) &= -\frac{(K^2+1)(24K^2-5)}{6K^2(4K^2-1)}.
    \end{align}
    \end{subequations}
    Given \eqref{KT2_exp_final}, the bifurcation equation expands as
    \begin{align}
\tilde{\mathcal{B}}(a,b,c_r) &= \begin{pmatrix} ac_r  + v_{3,0,0}(K)a^3 + v_{1,2,0}(K)ab^2 + \mathcal{O}\left(|a|+|b|+|c_r| \right)^4 \\ bc_r + w_{2,1,0}(K)a^2b +w_{0,3,0}(K)b^3  +  \mathcal{O}\left(|a|+|b|+|c_r| \right)^4 \end{pmatrix} = 0. \label{BtildeExpK}
\end{align}
We rescale variables as in \eqref{rescale3} and arrive at the new bifurcation equation
\begin{align}
\tilde{\tilde{\mathcal{B}}}(a,\tilde{b},\tilde{c_r}) = \begin{pmatrix} v_{3,0,0}(K) + \tilde{c_r} + v_{1,2,0}(K)\tilde{b}^2 + \mathcal{O}\left(|a|\left(1+|\tilde{b}|+|a|^2|\tilde{c_r}| \right)^4 \right) \\ w_{2,1,0}(K)\tilde{b} + \tilde{c_r}\tilde{b} + w_{0,3,0}(K)\tilde{b}^3 +  \mathcal{O}\left(|a|\left(1+|\tilde{b}|+|a|^2|\tilde{c_r}| \right)^4 \right)\end{pmatrix} = 0. \label{btildek}
\end{align}
When $a = 0$, we obtain a solution to the bifurcation equation provided
\begin{subequations}
\begin{align}
\tilde{c_r} = -v_{3,0,0}(K)-v_{1,2,0}(K)\tilde{b}^2,\quad\quad\quad& \label{ceqnK} \\ \big(w_{0,3,0}(K)-v_{1,2,0}(K)\big)\tilde{b}^3 - \big(v_{3,0,0}(K)-w_{2,1,0}(K)\big)\tilde{b} &= 0. \label{crbeqnsK}
\end{align}
\end{subequations}
Notably, the cubic for $\tilde{b}$ is simpler in the present case than in the $K = 3$ case. This is because  the leading-order Fourier modes of $u_r(x;a,b,c_r)$  in $\mathcal{T}_{\omega,2}$  have more resonant interactions with $\cos(x)$ and $\cos(Kx)$ when $K = 3$. The obvious solution of the cubic is $\tilde{b} = 0$. The other solutions satisfy
\begin{align}
\tilde{b}^2 = \frac{v_{3,0,0}(K)-w_{2,1,0}(K)}{w_{0,3,0}(K)-v_{1,2,0}(K)} = -\frac{K^2(4K^2-61)}{61K^2-4}.
\end{align}
For $K \geq 4$, the right-hand side of this equation is always negative. Thus, the remaining two solutions of the cubic are imaginary and will not contribute real-valued Wilton ripple solutions. 

If $\tilde{b} = 0$, necessarily $\tilde{c_r} = -v_{3,0,0}(K)$ by \eqref{ceqnK}. Then, we have 
\begin{align}
    \tilde{\tilde{\mathcal{B}}}(0,0,-v_{3,0,0}(K)) = 0,
\end{align}
by construction. The corresponding Jacobian is
\begin{align}
\frac{\partial \tilde{\tilde{\mathcal{B}}}}{\partial(\tilde{b},\tilde{c_{r}})}(0,0,-v_{3,0,0}(K)) = \begin{pmatrix} 0 & 1 \\ w_{2,1,0}(K)-v_{3,0,0}(K) & 0\end{pmatrix}, \label{jacK}
\end{align}
and a direct calculation confirms
\begin{align}
\textrm{det}\left(\begin{pmatrix} 0 & 1 \\ w_{2,1,0}(K)-v_{3,0,0}(K) & 0\end{pmatrix} \right) = \frac{(K^2+1)(4K^2-61)}{6(K^2-4)(4K^2-1)} > 0 \quad \textrm{for} \quad K \geq 4.
\end{align}
Thus, \eqref{jacK} is invertible, and the implicit function theorem gives us real analytic functions $\tilde{b}(a)$ and $\tilde{c_r}(a)$ such that $\tilde{\tilde{\mathcal{B}}}(a,\tilde{b}(a),\tilde{c_r}(a)) = 0$, $\tilde{b}(0) = 0$, and $\tilde{c_r}(0) = -v_{3,0,0}(K)$.

    For each $K \geq 4$, we have a unique $a$-parameter family of solutions to \eqref{kEqn_final} for sufficiently small $a \in \mathbb{R}$. These solutions travel at velocity
     \begin{align}
c(a) = c_0 + a^2\tilde{c_{r}}(a) = \frac{K^2}{K^2+1} + a^2\tilde{c_{r}}(a) , \label{cK}
     \end{align}
    and take the form
    \begin{align}
u(x;a) = a\Big(\cos(x) + \tilde{b}(a)\cos(Kx)\Big) + u_r\Big(x;a,a\tilde{b}(a),a^2\tilde{c_r}(a)\Big), \label{uK}
    \end{align}
     where $u_r \in \textrm{H}^4_{\textrm{per},\textrm{even}}(\mathbb{T})$ is real analytic in its parameters with Taylor expansion \eqref{ur_exp} and coefficients \eqref{u2_k}-\eqref{u02_k}. The functions $\tilde{b}(a)$ and $\tilde{c_r}(a)$ are real analytic for sufficiently small $a$ and satisfy
     \begin{align}
\tilde{b}(0) = 0 \quad \textrm{and} \quad  \tilde{c_r}(0) = \frac{(K^2+1)(5K^2-24)}{6K^2(K^2-4)}.
     \end{align}

    Unlike before, $\tilde{b}(0) = 0$ when $K \geq 4$. As a consequence, it is unclear whether $\tilde{b}(a)$ vanishes identically for sufficiently small $a$. If this is the case, then our Wilton ripple solutions degenerate into Stokes wave solutions with profiles that are orthogonal to $\cos(Kx)$. In the following lemma, we prove that this degeneracy is, in fact, not possible. An important step of our argument appeals to \cite{langer2024wilton}, where formal asymptotic computations of \eqref{cK}-\eqref{uK} were carried out to arbitrarily high order in $a$ via a Poincar\'{e}-Lindstedt method. 
\begin{lemma} For each $K \geq 4$, $\tilde{b}(a) \not\equiv 0$ for sufficiently small $a$.
\end{lemma}

\begin{proof}

    Fix $K \geq 4$. Since $\tilde{b}(a)$ and $\tilde{c}(a)$ are real analytic for sufficiently small $a$ and since \newline $u_r(x;a,a\tilde{b}(a),a^2\tilde{c}_r(a))$ is a composition of real analytic functions in its final three arguments, $c(a)$ and $u(x;a)$ are real analytic for sufficiently small $a$ with values in $\mathbb{R}$ and $\textrm{H}^4_{\textrm{per,even}}(\mathbb{T})$, respectively. Hence, for sufficiently small $a$, we have unique power series 
    \begin{subequations}
\begin{align}
c(a) &= \sum_{j = 0}^\infty c_j a^j, \quad \textrm{absolutely convergent in}~\mathbb{R}, \label{caexp}  \\
u(x;a) &=  \sum_{j = 1}^\infty u_j(x) a^j, \quad \textrm{absolutely convergent in H}^4_{\textrm{per,even}}(\mathbb{T}), \label{uaexp}
\end{align}
\end{subequations}
with $c_j \in \mathbb{R}$ and $u_j(x) \in\textrm{H} ^{4}_{\textrm{per,even}}(\mathbb{T})$. By \eqref{cK}, 
\begin{align}
c_0 = \frac{K^2}{K^2+1} \quad \textrm{and}  \quad c_2 = \frac{(K^2+1)(5K^2-24)}{6K^2(K^2-4)}.
\end{align}
Furthermore, since $\tilde{b}(a)$ is at least $\mathcal{O}\left(a\right)$ as $a \rightarrow 0$ and since $u_r(x;a,a\tilde{b}(a),a^2\tilde{c_r}(a))$ is at least\footnote{$u_r(x;a,b,c_r)$ is quadratic in $a$, $b$, and $c_r$ at leading order according to Case 2 in the Appendix.} $\mathcal{O}(a^2)$ as $a \rightarrow 0$, we conclude
\begin{align*}
u_1(x) = \cos(x).
\end{align*}
Explicit expressions of $c_j$ and $u_j(x)$ for larger values of $j$ require further expansions of the auxiliary and bifurcation equations.

Returning to the original traveling-wave equation \eqref{opEqn}, we have
\begin{align}
F(u(x;a),c(a)) = 0,   \quad \textrm{where} \label{ODE}
\end{align}
\begin{align}
    F(u,c) := (c+\mathcal{L})u + u^2, \quad \mathcal{L} := \partial_x^2 + \beta \partial_x^4, \quad \beta = 1/(1+K^2), \quad \textrm{and} \quad K \geq 4.
\end{align}
Since $c(a)$ and $u(x;a)$ are real analytic for sufficiently small $a$, $F(u(x;a),c(a))$ is also real analytic with values in $\textrm{L}^2_{\textrm{per,even}}(\mathbb{T})$. Thus, $F(u(x;a),c(a))$ admits a unique power series in $a$. Substituting expansions \eqref{caexp}-\eqref{uaexp} into \eqref{ODE} and using the Cauchy product to rearrange terms as needed\footnote{Rearrangement is permissible since $c(a)$ and $u(x;a)$ have absolutely convergent power series and since $F(u,c)$ is a real analytic map of $u$ and $c$ into $\textrm{L}^2_{\textrm{per,even}}(\mathbb{T})$.} yields the following power series of $F(u(x;a),c(a))$:
\begin{align}
F(u(x;a),c(a)) = \sum_{j = 1}^{\infty} \biggr( \sum_{i = 1}^{j}c_{j-i}u_i(x) + \mathcal{L}u_j(x) + \sum_{i = 1}^{j-1}u_{j-i}(x)u_i(x) \biggr)a^j, \label{Fseries}
\end{align}
which is absolutely convergent in $\textrm{L}^2_{\textrm{per,even}}(\mathbb{T})$ for sufficiently small $a$. From \eqref{ODE} and \eqref{Fseries}, we have
\begin{align}
\big(c_0 + \mathcal{L} \big)u_j(x) = -\sum_{i = 1}^{j-1}\Big(c_{j-i} + u_{j-i}(x) \Big)u_i(x) \quad \textrm{for} \quad j \in \mathbb{N}, \label{ODEHeir}
\end{align}
 where, for notational convenience, the sum on the right-hand side collapses when $j = 1$. 

In addition to \eqref{ODE}, we also have the following for $\omega \in \{1,K\}$:
\begin{align}
\mathcal{F}_\omega\Big[F(u(x;a),c(a))\Big] &= 0, \nonumber \\
\mathcal{F}_\omega\left[c(a)u(x;a) + u(x;a)^2\right] + \mathcal{F}_\omega\Big[\mathcal{L} u(x;a) \Big] &= 0, \nonumber \\
\mathcal{F}_\omega\Big[c(a)u(x;a) + u(x;a)^2 \Big] - c_0\mathcal{F}_\omega\Big[u(x;a)\Big] &= 0, \nonumber \\
\mathcal{F}_\omega\Big[\Big(c(a)-c_0+u(x;a) \Big) u(x;a)\Big] &= 0, \label{FODE}
\end{align}
where the third equality follows from the self-adjointness of $\mathcal{L}$ and the definition of $c_0$. The left-hand side of \eqref{FODE} is real analytic for sufficiently small $a$ with values in $\mathbb{R}$.  Substituting expansions \eqref{caexp}-\eqref{uaexp} into \eqref{FODE} yields the unique power series
\begin{align}
\sum_{j = 2}^{\infty} \mathcal{F}_{\omega}\Big[ \sum_{i = 1}^{j-1} \big(c_{j-i} + u_{j-i}(x) \big)u_i(x)\Big]a^j = 0,
\end{align}
which converges absolutely in $\mathbb{R}$ for sufficiently small $a$. It follows that
\begin{align}
\mathcal{F}_\omega\Big[\sum_{i = 1}^{j-1} \big(c_{j-i} + u_{j-i}(x) \big)u_i(x) \Big] = 0 \quad \textrm{for} \quad j \in \mathbb{N} \setminus \{1\} \quad \textrm{and} \quad \omega \in \{1,K\}. \label{ODESolv}
\end{align}

Together, \eqref{ODEHeir} and \eqref{ODESolv} are (up to notation) exactly the hierarchy of equations studied in the formal work of \cite{langer2024wilton} to obtain arbitrarily high-order asymptotic expansions of $c(a)$ and $u(x;a)$ via the Poincar\'{e}-Lindstedt method. We have hence demonstrated that these formal computations are consistent with the functional-analytic framework developed in this paper. That is, the formal, high-order asymptotic expansions computed by \cite{langer2024wilton} converge for sufficiently small $a$ to the $c(a)$ and $u(x;a)$ obtained in this work. Having established the rigorous validity of these expansions, we turn to \cite[Lemma 5]{langer2024wilton}, which shows $\mathcal{F}_K[u(x;a)] \not= 0$ and, more specifically, $\mathcal{F}_K[u(x;a)] = \mathcal{O}\left(a^{K-2} \right)$ as $a \rightarrow 0$. Consequently, $\tilde{b}(a) \not\equiv 0$ for sufficiently small $a$ and, in fact, $\tilde{b}(a) = \mathcal{O}\left(a^{K-3}\right)$ as $a \rightarrow 0$.
\end{proof}

With Lemma 1 in hand, we conclude for each $K \geq 4$ that \eqref{uK} represents a non-degenerate Wilton ripple solution  depending real analytically on $a$ with values in $\textrm{H}^4_{\textrm{per,even}}(\mathbb{T})$. As in the previous cases, we apply a bootstrap argument to \eqref{kEqn_final} to promote $u(x;a)$ to the more regular space $\textrm{H}^{\infty}_{\textrm{per,even}}(\mathbb{T})$, completing the proof of Theorem 1.
   
    \section*{Conclusion}

In this paper, we have rigorously established the existence of all $1\!:\!K$ Wilton ripple solutions of the Kawahara equation for each $K \in \mathbb{N} \setminus \{1\}$. These solutions are defined as small-amplitude, periodic, traveling-wave solutions that bifurcate from a linear combination of $\cos(x)$ and $\cos(Kx)$ in the limit of vanishing amplitude, where the Bond parameter $\beta = 1/(1 + K^2)$ induces a degeneracy in the linearized dispersion relation. The resulting bifurcation is of codimension one and leads to a two-dimensional kernel, making the construction of these solutions substantially more delicate than for standard Stokes waves.

Our proof proceeds via a Lyapunov--Schmidt reduction, yielding a reduced bifurcation equation with real analytic dependence on a small-amplitude parameter. The central technical challenge is to verify that the coefficient function $b(a)$ in the Wilton ripple expansion~\eqref{wilton_exp} does not vanish identically. This requires a precise understanding of the asymptotic structure of the Wilton ripple solutions and their dependence on $K$. By rigorously justifying the high-order formal expansions developed in~\cite{langer2024wilton}, we establish the nontriviality of $b(a)$ and thereby conclude the existence of genuine $1:K$ Wilton ripple branches for each $K$.

Previous rigorous existence results for Wilton ripples have largely been restricted to the 1:2 resonance case \cite{akers2021wiltonB,akers2021wiltonA, reeder1981wiltona,reeder1981wiltonb} or relied on bifurcation frameworks that allow equation parameters such as the Bond number $\beta$, vorticity, or stratification to vary \cite{ehrnstrom2015trimodal,ehrnstrom2019bifurcation,kozlov2018modal,martin2013existence,wang2025existence}. While these approaches yield families of multimodal solutions organized in parameter sheets, they do not directly imply the existence of Wilton ripples in the classical sense, where equation parameters are fixed \cite{akers2025ripples}. The results presented here bridge that gap, confirming the existence of discrete Wilton branches at fixed $\beta$ for all $K$, in a setting that admits full analytic control.

The framework developed in this work is expected to generalize to other nonlinear dispersive models, including the full gravity–capillary water wave equations. A major technical hurdle, however, is that, for each $K \in \mathbb{N} \setminus \{1\}$, the Wilton ripple solutions must be expanded to increasingly high orders in the small-amplitude parameter in order to detect the resonant $\cos(Kx)$ mode. While such computations remain out of reach in the full water wave setting, the formal expansions carried out in~\cite{langer2024wilton} for weakly nonlinear models raise the intriguing possibility that analogous expansions could eventually be constructed for more physically realistic systems. Should that become possible, the method of proof outlined here may offer the first rigorous existence theory for all $1\!:\!K$ Wilton ripples in the full gravity–capillary water wave equations. Future work may also explore whether these fixed-parameter branches of Wilton ripples lie embedded within the bifurcation sheets arising in more recent treatments of multimodal solutions.

     \section{Appendix: The Taylor Coefficients $u_r^{(i,j,k)}(x)$}
\subsection{Case 1: $K = 2$} 
We begin with the expansion of $u_r(x;a,b,c_r)$ when $K = 2$. Substituting \eqref{ur_exp} into \eqref{Aeqn} and equating collective powers in $a$, $b$, and $c_r$ yields a hierarchy of inhomogeneous differential equations that uniquely determines $u_r^{(i,j,k)}(x)$ for increasingly larger $i$, $j$, and $k$. Calculating these equations is straightforward but tedious, so we appeal to Mathematica for help. The interested reader can consult \emph{CompanionToKawaharaWiltonRipples.nb} for further details. 

The differential equations for $u_r^{(i,j,k)}(x)$ coincide when $i +j + k = 1$  and are given by 
\begin{align}
(c_0+\mathcal{L})u^{(i,j,k)}(x) = 0.
\end{align}
As $u^{(i,j,k)}(x) \in \textrm{Null}\{c_0+\mathcal{L}\}^\perp$, we conclude
\begin{align}
u^{(i,j,k)}(x) = 0 \quad \textrm{for} \quad i+j+k =1. \label{urlin1}
\end{align}
In other words, there are no first-order terms in the Taylor Series of $u_r(x;a,b,c_r)$.

When $i + j + k = 2$, we obtain six differential equations. In particular, 
\begin{subequations}
\begin{align}
(c_0+\mathcal{L})u^{(2,0,0)}(x) &= -\frac12, \\
(c_0+\mathcal{L})u^{(1,1,0)}(x) &= -\cos(3x), \\
(c_0+\mathcal{L})u^{(0,2,0)}(x) &= -\frac12-\frac12\cos(4x), \\
(c_0+\mathcal{L})u^{(0,1,1)}(x) &= 0, \\
(c_0+\mathcal{L})u^{(0,0,2)}(x) &= 0, \\
(c_0+\mathcal{L})u^{(1,0,1)}(x) &= 0.
\end{align}
\end{subequations}
By undetermined coefficients (or an equivalent method), one finds
\begin{subequations}
\begin{align}
u^{(2,0,0)}(x) &= -\frac{1}{2c_0} , \label{2u200}\\
u^{(1,1,0)}(x) &= -\frac{1}{c_0 - 9 + 81\beta}\cos(3x), \label{2u110} \\
u^{(0,2,0)}(x) &=  -\frac{1}{2c_0} - \frac{1}{2\left(c_0 - 16 +256\beta \right)}\cos(4x). \label{2u020}
\end{align}
\end{subequations}
The remaining coefficient functions are zero by the same logic as in the $i+j+k=1$ case. 
\subsection{Case 2: $K \geq 3$} 
We proceed as in Case 1 but leave the details to Mathematica. When $i + j + k = 1$, we find
\begin{align}
u^{(i,j,k)}(x) = 0. \label{urlin2}
\end{align}
When $i + j + k = 2$, we have instead
\begin{subequations}
\begin{align}
u^{(2,0,0)}(x) &= -\frac{1}{2c_0} - \frac{1}{2\left(c_0 - 4 + 16\beta \right)}\cos(2x) \label{u2_k}\\
u^{(1,1,0)}(x) &= -\frac{1}{c_0 - (K-1)^2 + \beta(K-1)^4}\cos((K-1)x) - \frac{1}{c_0 - (K+1)^2 + \beta(K+1)^4}\cos((K+1)x), \label{u11_k}\\
u^{(0,2,0)}(x) &= -\frac{1}{2c_0} - \frac{1}{2\left(c_0 - 4K^2 + 16\beta K^4 \right)}\cos(2Kx), \label{u02_k}
\end{align}
\end{subequations}
and the remaining coefficient functions are zero. \\\\
\emph{Remark}: If $K = 2$, the coefficients of $\cos(2x)$ in \eqref{u2_k} and of $\cos((K-1)x)$ in \eqref{u11_k} are singular. This is to be expected, as $u^{(i,j,k)}(x)$ must be orthogonal to $\textrm{Null}\{c_0+\mathcal{L}\} = \textrm{Span}\{\cos(x),\cos(2x)\}$ when $K = 2$. It is for this reason that this case is considered separately. 



\begin{thebibliography}{999}

\bibitem{akers2025ripples}
B. F. Akers. “On ripples: bifurcations of resonant bimodal traveling waves”. In: \emph{Nonlinearity} 38.6 (2025), pp. 065005.

\bibitem{akers2021wiltonB}
B. Akers and D. P. Nicholls. “Wilton ripples in weakly nonlinear dispersive models of water
waves: Existence and analyticity of solution branches”. In: \emph{Water Waves} 3.1 (2021), pp. 25–47.

\bibitem{akers2021wiltonA}
B. Akers and D. P. Nicholls. “Wilton Ripples in weakly nonlinear models of water waves:
Existence and computation”. In: \emph{Water Waves} 3.3 (2021), pp. 491–511.

\bibitem{alyousef2022new}
H. A. Alyousef et al. “New Periodic and Localized Traveling Wave Solutions to a Kawahara-
Type Equation: Applications to Plasma Physics”. In: \emph{Complexity} 2022.1 (2022), p. 9942267.

\bibitem{assas2009new}
L. M. B. Assas. “New exact solutions for the Kawahara equation using Exp-function method”.
In: \emph{Journal of Computational and Applied Mathematics} 233.2 (2009), pp. 97–102.

\bibitem{chow2012methods}
S.-N. Chow and J. K. Hale. \emph{Methods of Bifurcation Theory}. Vol. 251. Springer Science \& Business Media, 2012.

\bibitem{coleman2012calculus}
R. Coleman. \emph{Calculus on Normed Vector Spaces}. Springer, 2012.

\bibitem{creedon2021high}
R. P. Creedon, B. Deconinck, and O. Trichtchenko. “High-frequency instabilities of the Kawahara equation: a perturbative approach”. In: \emph{SIAM Journal on Applied Dynamical Systems}
20.3 (2021), pp. 1571–1595.

\bibitem{ehrnstrom2015trimodal}
M. Ehrnstr{\"o}m and E. Wahlen. ``Trimodal steady water waves''. In: \emph{Archive for Rational Mechanics and Analysis} 216.2 (2015), pp. 449-471.

\bibitem{ehrnstrom2019bifurcation}
M. Ehrnstr\"{o}m et al. “On the bifurcation diagram of the capillary-gravity Whitham equation”.
In: \emph{Water Waves} 1.2 (2019), pp. 275–313.

\bibitem{haragus2006spectral}
M. Haragus, E. Lombardi, and A. Scheel. “Spectral stability of wave trains in the Kawahara
equation”. In: \emph{Journal of Mathematical Fluid Mechanics} 8 (2006), pp. 482–509.

\bibitem{haupt1988modeling}
S. E. Haupt and J. P. Boyd. “Modeling nonlinear resonance: A modification to the stokes’
perturbation expansion”. In: \emph{Wave Motion} 10.1 (1988), pp. 83–98.

\bibitem{henderson1987experiments}
D. M. Henderson and J. L. Hammack. “Experiments on ripple instabilities. Part 1. Resonant
triads”. In: \emph{Journal of Fluid Mechanics} 184 (1987), pp. 15–41.

\bibitem{jones1989small}
M. C. W. Jones. “Small amplitude capillary-gravity waves in a channel of finite depth”. In:
\emph{Glasgow Mathematical Journal} 31.2 (1989), pp. 142–160.

\bibitem{kawahara1972}
T. Kawahara. “Oscillatory solitary waves in dispersive media”. In: \emph{J. Phys. Soc. Jpn.} 33 (1972), pp. 1015–1024.

\bibitem{kozlov2018modal}
V. Kozlov and E. Lokharu. ``N-modal steady water waves with vorticity''. In: \emph{Journal of Mathematical Fluid Mechanics} 20.2 (2018), pp. 853-867. 

\bibitem{langer2024wilton}
R. Langer, O. Trichtchenko, and B. Akers. “Wilton Ripples with High-Order Resonances in Weakly Nonlinear Models”. In: \emph{Water Waves} (2024), pp. 1–30.

\bibitem{mancas2019traveling}
S. C. Mancas. “Traveling wave solutions to Kawahara and related equations”. In: \emph{Differential
Equations and Dynamical Systems} 27.1 (2019), pp. 19–37.

\bibitem{martin2013existence}
C. I. Martin and B.-V. Matioc. “Existence of Wilton ripples for water waves with constant
vorticity and capillary effects”. In: \emph{SIAM J. Appl. Math.} 73.4 (2013), pp. 1582–1595.

\bibitem{maspero2024full}
A. Maspero and A. M. Radakovic. “Full description of Benjamin-Feir instability for generalized
Korteweg-de Vries equations”. In: \emph{arXiv:2404.06172} (2024).

\bibitem{mcgoldrick1970experiment}
L. F. McGoldrick. “An experiment on second-order capillary gravity resonant wave interac-
tions”. In: \emph{Journal of Fluid Mechanics} 40.2 (1970), pp. 251–271.

\bibitem{reeder1981wiltona}
J. Reeder and M. Shinbrot. “On Wilton ripples, I: Formal derivation of the phenomenon”. In:
\emph{Wave Motion} 3.2 (1981), pp. 115–135.

\bibitem{reeder1981wiltonb}
J. Reeder and M. Shinbrot. “On Wilton ripples, II: rigorous results”. In: \emph{Archive for Rational
Mechanics and Analysis} 77.4 (1981), pp. 321–347.

\bibitem{sirendaoreji2004new}
Sirendaoreji. “New exact travelling wave solutions for the Kawahara and modified Kawahara
equations”. In: \emph{Chaos, Solitons and Fractals} 19.1 (2004), pp. 147–150.

\bibitem{stokes1847theory}
G. G. Stokes. “On the theory of oscillatory waves”. In: \emph{Trans. Cam. Philos. Soc.} 8 (1847), pp. 441–455.

\bibitem{toland1985bifurcation}
J. F. Toland and M. C. W. Jones. ``The bifurcation and secondary bifurcation of capillary-gravity waves". In: \emph{Proceedings of the Royal Society of London. A. Mathematical and Physical Sciences} 399.1817 (1985), pp. 391-417.  

\bibitem{trichtchenko2018stability}
O. Trichtchenko, B. Deconinck, and R. Koll\'{a}r. “Stability of periodic traveling wave solutions
to the Kawahara equation”. In: \emph{SIAM Journal on Applied Dynamical Systems} 17.4 (2018),
pp. 2761–2783.

\bibitem{trichtchenko2016instability}
O. Trichtchenko, B. Deconinck, and J. Wilkening. ``The instability of Wilton ripples''. \emph{Wave Motion} 66 (2016), pp. 146-155. 

\bibitem{wang2025existence}
J. Wang, F. Xu, and Y. Zhang. ``The existence of stratified linearly steady two-mode water waves with stagnation points''. In: \emph{Journal of Mathematical Fluid Mechanics} 27.1 (2025), pp. 13. 

\bibitem{wilton1915lxxii}
J. R. Wilton. “LXXII. On Ripples”. In: \emph{The London, Edinburgh, and Dublin Philosophical
Magazine and Journal of Science} 29.173 (1915), pp. 688–700.

\bibitem{ye2022all}
F. Ye et al. “All traveling wave exact solutions of the Kawahara equation using the complex
method”. In: \emph{Axioms} 11.7 (2022), p. 330.

\bibitem{zakaria2004wilton}
K. Zakaria. “Wilton ripples between two uniform streaming magnetic fluids”. In: \emph{International
Journal of Non-Linear Mechanics} 39.6 (2004), pp. 1051–1066.
\end{thebibliography}
\end{document}